\documentclass[a4paper,11pt,english]{amsart}
\usepackage{amsthm,amsmath,amssymb}
\usepackage[T1]{fontenc}
\usepackage[utf8]{inputenx}
\usepackage{graphicx}
\usepackage{geometry}
\usepackage{tikz}
\usepackage{hyperref}
\usepackage{leftidx}
\usepackage{xcolor}
\usepackage{pifont}
\usepackage{babel}
\usepackage{enumerate}
\usepackage{soul} 
\usepackage{comment}

\newtheorem{theorem}{Theorem}

\newtheorem{lemma}[theorem]{Lemma}
\newtheorem{proposition}[theorem]{Proposition}
\newtheorem{definition}{Definition}

\newtheorem{remark}{Remark}

\newtheorem{corollary}[theorem]{\bf Corollary}
 
\usepackage{thmtools}
\usepackage{cleveref}

\newcommand{\E}{\mathbb{E}}
\newcommand{\N}{\mathbb{N}}

\renewcommand{\P}{\mathbb{P}}
\newcommand{\R}{\mathbb{R}}
\newcommand{\Z}{\mathbb{Z}}

\newcommand{\Pbarre}{\overline{\mathbb{P}}}

\newcommand{\1}{1\hspace{-1.3mm}1}

\title{An asymptotic shape theorem for additive random linear growth models}
\author{Aurelia Deshayes  \and Pierrick Siest}
\date{\today}


\begin{document}
\begin{abstract}
In this paper, we define a class of additive random growth models whose growth is at least and at most linear and prove an asymptotic shape theorem for these models. This proof generalizes already known proofs for the classical contact process (\cite{harris74}, \cite{DG82}) or some of its variants (contact process on supercritical random environment\cite{GMPCRE} or contact process with aging~\cite{CPA}) and allows us to obtain conjectured asymptotic shape theorems for Richardson's model with stirring and the contact process with stirring~\cite{MMS}.
\end{abstract}
\maketitle

\section{Introduction}

In 1974, Harris~\cite{harris74} introduced the \emph{classical contact process} as an interacting particle system modeling the spread of an epidemic through the grid $\Z^d$ with the following dynamics: an infected site recovers at rate 1, and infects its neighbors at rate $\lambda$. Harris showed that, for any dimension $d$, this process exhibits a phase transition. This means that there exists a non-trivial value $\lambda_c(d)$, called critical value, such that if $\lambda>\lambda_c(d)$ (respectively $\lambda<\lambda_c(d)$), then the probability that, starting from a single infected point, the infection persists over time is positive (respectively zero). In dimension~1, thanks to a coupling with an oriented percolation, Durrett and Griffeath~\cite{DG83} proved that, in the supercritical regime, the rightmost infected particle has linear speed. In higher dimensions, the analogous result is called \emph{asymptotic shape theorem}: denoting by $H_t$ the set of particles infected at least once before time $t$, is there a deterministic set $B$ such that $\frac{H_t}{t}$ converges to $B$ (in a sense to be made precise)? Durrett and Griffeath~\cite{DG82} proved this result for large infection rates $\lambda$ only, since they needed the linear growth of the contact process, which at the time were only proved for large $\lambda$. Then, Bezuidenhout-Grimmett~\cite{BG90} made a smart block construction and proved these estimates to be true in the whole supercritical regime, thus completing the proof of the asymptotic shape theorem for the classical contact process. An overview can be found in~\cite{Dur91} or~\cite{lig99}. 

Since then, numerous extensions have appeared in the literature: the \emph{two stage contact process} by Krone~\cite{krone}, the \emph{boundary modified contact process} by Durrett and Schinazi~\cite{durschi}, the \emph{contact process in randomly evolving environment} by Broman~\cite{broman}, Remenik~\cite{remenik}, Steif and Warfheimer~\cite{steifwar}, the \emph{contact process with aging}~\cite{CPA}, the \emph{contact process with stirring}~\cite{MMS}, \textit{etc...} For all these processes, we can exhibit a property of interest $P$ about the state of a site of $\Z^d$, then define $H_t$ as the subset of $\Z^d$ consisting of all points that have satisfied the property at least once before time $t$. All these processes are linear random growth models in the following sense: there are deterministic non-empty compact sets $H_-$ and $H_+$ such that, for $t$ large enough, with high probability:
\[tH_- \subset H_t \subset tH_+.\]
The first inclusion expresses the \emph{at least linear growth} and the second one the \emph{at most linear growth}.
These models also share characteristics such as attractivity, phase transition phenomenon \textit{etc}. In this paper, the goal is to highlight their common properties which allow to prove that a general class of linear random growth models satisfies an asymptotic shape theorem. We call \emph{hitting time} of $x$ the first time $t(x)$ that $x$ satisfies the property $P$. Models for which once a site $x$ satisfies the property $P$, so it does forever, are called \emph{permanent} models. It includes, for example, first passage percolation, introduced by Hammersley and Welsh~\cite{hammwelsh}, considering the property of "being wet". For permanent models, Kingman's theorem is generally used on the sequence of hitting times $(t(nx))_{n\in \N}$ for all $x\in \Z^d$, thanks to the hypotheses of sub-additivity, stationarity and integrability. In our case of random linear growth models, we consider non-permanent models: the fact that the property may disappear means that the standard integrability conditions are not satisfied, since $t(nx)$ can be infinite with positive probability. If we condition on the property's survival, then stationarity and subadditivity properties can be lost. To overcome such lacks, we introduce a suitable quantity $\sigma(x)$ called \emph{essential hitting time} of $x$, inspired by~\cite{GMPCRE}. It can been seen as a regeneration time: it is a time when the site $x$ satisfies the property, and spreads it forever. It turns out that this function $\sigma$ satisfies, under some hypotheses,  adequate stationarity properties as well as the almost-subaddivity conditions involved in Kesten and Hammersley's theorem (see \cite{kestentrick} and \cite{hamm74}), a well-known extension of Kingman's seminal result. Therefore, we can use the same techniques as the ones involved in~\cite{CPA}, in order to prove an asymptotic shape theorem, that is the existence of a norm $\mu$ such that, for all $\varepsilon>0$, almost surely on the event "the property persists over time", for all $t$ large enough, we have that
\[(1-\varepsilon)B_{\mu}\subset\frac{\tilde{H}_t}t\subset(1+\varepsilon)B_{\mu},\]
with $\tilde{H}_t=H_t+[0,1)^d$, and $B_\mu$ the unit ball associated to $\mu$.

Our techniques work with a general state space, provided that the process is well defined, constructible from Poisson processes, and satisfies the growth assumptions. However, they require that the process be additive. This assumption is discussed at the end of Section 2. It excludes many examples of the application of our theorem, and we believe that it can be removed.

In Section~\ref{section_class}, we define our class of random linear growth models, and give examples of such models. We then state the general asymptotic shape theorem. We deduce from it asymptotic shape theorems that were conjectured for 
Richardson's model with stirring and contact process with stirring~\cite{MMS}. 
The proof of the general asymptotic shape theorem is cut into two parts; in Section~\ref{section_essential}, we define the essential hitting time, and state its properties; in Section~\ref{section_thm}, we prove the general asymptotic shape theorem for the essential hitting time, and deduce from it the general asymptotic shape theorem for the hitting time. 

\section{A class of random linear growth models}\label{section_class}
\subsection{Notation and definitions} 
\subsubsection{State space and property of interest} We work on interacting particle systems on $\Z^d$. For $x$ and $y$ in $\Z^d$, we say that $x$ and $y$ are neighbors, and we denote by $x\sim y$, if $||x-y||_1=\sum_{i=1}^d |x_i-y_i|=1$. Let $S$ be a countable well-ordered set, that is, a totally ordered set with the property that every subset of $S$ has a smallest element. We call $S$ the \emph{states space}: it contains the possible states (or types) of a particle $x\in\Z^d$.
We endow $S$ with the discrete topology, and we denote by $S^{\Z^d}$ the set of mappings $\xi:\Z^d\to S$. The set $S^{\Z^d}$ is considered as the space of configurations:
\[\xi\in S^{\Z^d} \Leftrightarrow \left\{\begin{array}{cccc}
                                          \xi: &\Z^d &\to &S \\
                                         & x &\to &\xi(x)
                                         \end{array}.
                                        \right.
                                        \] 
We define a \emph{property of interest} $P$ as a subset of $S$ which satisfies the following properties:
\begin{itemize}
    \item $P$ is a non-trivial subset of $S$ (that is, $P\notin \{\emptyset,S\}$).
    \item  $P$ is an \emph{upper set}, that is, if $s\in P$ and $s\le s'$, then $s'\in P$.
\end{itemize}
We can think of the property $P$ as 'to be in the maximal state' when $S$ is bounded of 'to be in a state bigger than some $p_0$' (for example 'to be infected' for the classical contact process). By definitions of $S$ and $P$, we can define, for all $x\in\Z^d$, the smallest configuration $\Delta_x$ in which only site $x$ satisfies the property $P$: for all $y\in \Z^d$,
\begin{align} \label{Définition configuration minimale}
\begin{cases}
\Delta_x (y)= \min(S) &\text{if}~ y\neq x \\
\Delta_x (x)= \min(P) 
\end{cases}
\end{align}

In the same way, for $x\in\Z^d$ and $p\in P$ we define $\Delta_x^p$ such that $\Delta_x^p(y)=\min(S)$ for $y\neq x$ and $\Delta_x^p(x)=p$.
Note that, since $P\ne S$ and that $P$ is an upper set, we have $\min(S)\notin P$.
Finally, we say that a configuration $\xi$ has a \emph{finite support} if the set $\{x\in \Z^d: \xi(x)\in P\}$ is finite.

\subsubsection{Markov process} 

We work with a continuous time space $\R_+$. 
Let $(\xi_t)_{t\in \R_+}$ be a Feller Markov process taking values in $S^{\Z^d}$. This process gives the state of each site at each time. We denote by $(\mathcal{F}_t)_{t\in \R_+}$ its associated filtration. For a configuration $f\in S^{\Z^d}$, $(\xi_t^f)$ is the process starting from initial configuration $\xi_0=f$. For $x\in \Z^d$, $(\xi_t^x)$ is the process $(\xi_t^{\Delta_x})$, starting from $\Delta_x$ the minimal configuration previously defined. If the initial configuration is not specified, it means that we start with the configuration $\Delta_{0_{\Z^d}}$. 

We define the processes $(A_t)_{t\in \R_+}$ and $(H_t)_{t\in \R_+}$, taking values in $\mathcal P (\Z^d)$, by
\[A_t=\{y\in\Z^d,~\xi_t(y)\in P\}\quad\text{and}\quad H_t=\bigcup_{s\le t} A_s;\]
$A_t$ is the set of points which satisfies the property $P$ at time $t$, which has no reason to be non-decreasing, and $H_t$ is the set of points which have satisfied $P$ before time $t$. The process $(H_t)$ is non-decreasing in the sense that once $x$ belongs to $H_{t_0}$, it belongs to $H_t$ for all $t>t_0$. The processes $(A_t^f),(H_t^f),(A_t^y),(H_t^y)$ are defined as before, according to the initial configuration.

In general, the process $(A_t)$ gives only partial information: if we want to know the set $A_{t+dt}$ of sites satisfying the property at time $t+dt$, we need the full knowledge of $\xi_t$. Therefore, we should be careful when using the Markov property with the process $(A_t)$. 

\subsubsection{Interacting particle system}

We are interested in a subclass of \emph{nearest neighbor interacting (NNI)} particle systems; a general definition and the proof of their existence in the case where $S$ is compact can be found in the book~\cite{IPS}. Here the set $S$ may not be compact; we do not prove that such interacting particle systems exist in general, it has to be done for each particular case (for example, it is done in~\cite{CPA} where $S=\N$).

Let us introduce some notation. Let $\mathcal{C}(S^{\Z^d})$ be the set of continuous functions from $S^{\Z^d}$ to $\R$, and $\mathcal{C}_0(S^{\Z^d})$ be the subset of continuous functions depending only on finitely many coordinates of $\xi$. For $x,y\in \Z^d$, $s\in S$, $\bar s=(s_1,s_2)\in S^{2}$ and $\xi \in S^{\Z^d}$, we set: 
\begin{itemize}
    \item $\xi_{(x,s)}$, the configuration identical to $\xi$, except for $\xi(x)$ which is replaced by $s$.
    \item $\xi_{(x,y,\bar s)}$, the configuration identical to $\xi$, except for $\xi(x)$ (resp. $\xi(y)$) which is replaced by $s_2$ (resp. $s_1$).
\end{itemize}
We study Feller Markov processes on $S^{\Z^d}$, which we suppose to be well-defined and which have a generator of the form:
\begin{align} \label{Forme du générateur}
    \mathcal{A}f(\xi) &=\sum_{x\in\Z^d}\sum_{s\in S} c(x,\xi,s) \left(f(\xi_{(x,s)})-f(\xi)\right)\\
    &+\sum_{\substack{{x,y\in\Z^d}\\ x\sim y}}\sum_{\bar s\in S^2} c(x,y,\xi,\bar s) \left(f(\xi_{(x,y,\bar s)})-f(\xi)\right),\notag
\end{align}
for $f\in \mathcal{C}_0(S^{\Z^d})$ and $\xi\in S^{\Z^d}$. The quantity $c(x,\xi,s)\ge 0$ (resp. $c(x,y,\xi,\bar s)\ge 0$) is the intensity of a jump $\xi \to\xi_{(x,s)}$ (resp. $\xi \to\xi_{( x,y,\bar s)}$), and depends on $\xi$ only through $\{\xi(z), z\sim x\}$ and $\xi(x)$  (resp. through $\{\xi(z), z\sim x\text{ or } z \sim y\}$). 
For simplicity, we take $c(x,\xi,s)=0$ if $\xi(x)=s$ and $c(x,y,\xi,\bar s)=0$ if $\xi(x)=s_2$ or $\xi(y)=s_1$. We also suppose that the transition rates are \emph{symmetric}, that is:
\begin{align}
    c(x,\xi,s)=c(x,\xi(-\ \cdot),s)\quad \text{and}\quad c(x,y,\xi,\bar s)=c(0,y-x,\xi(-\ \cdot),\bar s). \label{symmetry property}
\end{align}

Note that in most of the models considered, only one site changes its state at a time, so the second term of the generator is zero. This term appears for models which incorporate stirring \cite{MMS}. 

We will suppose that our interacting particle systems are \emph{translation invariant} in the sense that for all $x,y\in \Z^d$, $\xi\in S^{\Z^d}$, $s\in S$ and $\bar s\in S^2$, we have 
\begin{align}
\label{Invariance par translation}
    c(x,\xi,s)=c(0,\xi(.-x),s) \quad\text{and}\quad c(x,y,\xi,\bar s)=c(0,y-x,\xi(.-x),\bar s).
\end{align}

\subsubsection{Construction with Poisson point processes} \label{Subsubsection Construction with Poisson processes}

We say that the process $(\xi_t)$ \emph{can be constructed from (marked) Poisson point processes} if there are families $(\mathcal N_{e})_{e\in \E^{d}}$ (where $\E^d$ is the set of oriented edges between nearest neighbors), and $(\mathcal N_{x})_{x\in \Z^{d}}$ of mutually independent Poisson point processes on $\R^+\times \mathcal{M}$, where $\mathcal{M}$ is a finite set of marks, such that:
\begin{itemize}
    \item the set of jump times of the process $(\xi_t)$ is a subset of the times given by these families,
    \item if $(T,M)\in \mathcal N_{e}$ (resp. $(T,M)\in \mathcal{N}_{x}$) is a jump time of $(\xi_t)$ of mark $M\in\mathcal{M}$, then we have $\xi_T=f_{e}(\xi_{T^-},M)$ (resp. $x_T=f_x(\xi_{T^-},M)$), where $f_{e}$ and $f_x$ are deterministic functions.
\end{itemize}
For many models, such a construction is done to obtain a graphical representation, as introduced by Harris in the case of $S=\{0,1\}$. Foxall~\cite{Foxall2016} constructed this for any additive model with finite state space. To our knowledge, there is no general result for the case $|S|=\infty$ and it has to be done in all of such cases. See~\cite{CPA} for an example of such a construction in the case where $S=\N$.

\begin{remark} Conversely, we could have built a Feller Markov process on $S^{\Z^d}$ from a family of Poisson point processes like it is done by Swart (\cite{swart}, Section $4.3$). When possible, he constructs the process from a Poisson point process $\omega$ on $\mathcal{G}\times \R_+$, where $\mathcal{G}$ is a set of local maps from $S^{\Z^d}$ to itself such that the generator $\mathcal{A}$ of $(\xi_t)$ can be written as:
$$\forall \xi\in S^{\Z^d},~\mathcal{A}f(\xi)=\sum_{m\in\mathcal{G}} r_m\left[f(m(\xi))-f(\xi)\right],$$
where for all $m\in \mathcal G$, $r_m$ is a positive constant.
\end{remark}


\subsubsection{Additive models}
The classical contact process and some of its variations are \emph{additive} processes, that is:  for all $f,g\in S^{\Z^d}$, the processes $(\xi_t^{f})$, $(\xi_t^{g})$ and $(\xi_t^{f\vee g})$ can be coupled (using the same family of Poisson point processes) in such a way that
\begin{align}
\forall t\in \R_+,\quad \xi_t^{f}\vee \xi_t^{g}=\xi_t^{f\vee g}, \label{Propriété d'additivité}
\end{align}
where for all $x\in \Z^d$, $(f\vee g)(x)=\max[f(x),g(x)]$. Harris~\cite{harris78} proved that additive processes can be constructed graphically based on Poisson point processes.

Let us summarize the class $\mathcal{C}$ of models we are considering:

\begin{definition} \label{classeC} Let $(\xi_t)$ be a Feller Markov process taking values in $S^{\Z^d}$, and $P$ a property of interest. We say that $((\xi_t),P)$ belongs to the class $\mathcal C$ if it satisfies the following conditions:
\begin{itemize} 
 \item The generator of the process $(\xi_t)$ is of the form~\eqref{Forme du générateur}, and verifies \eqref{Invariance par translation} (translation invariance). 
 \item The process $(\xi_t)$ can be constructed with Poisson point processes such that $(\xi_t)$ is additive.
 \item If there exists $t_0\in \R_+$ such that $A_{t_0}=\emptyset$, then for all $t\ge t_0$, we have $A_t=\emptyset$.
\end{itemize}
\end{definition}

We state two properties that are verified by models of the class $\mathcal{C}$.

\begin{lemma} 
Let $((\xi_t),P)\in \mathcal{C}$.
    \begin{itemize}
        \item The process $(\xi_t)$ is \emph{attractive}, i.e. for all $f,g\in S^{\Z^d}$, the processes $(\xi_t^f)$ and $(\xi_t^g)$, can be constructed with the same family of Poisson processes such that :
            \begin{align}
            \text{if }f\le g,\text{ then for all }t\ge 0,\quad \xi_t^f\le \xi_t^g. \label{Propriété d'attractivité}
            \end{align}
        \item The configurations $\Delta_x$, $x\in \Z^d$, defined in \eqref{Définition configuration minimale}, are the worst initial configurations for the spread of the property $P$, in the sense that for all $f \in S^{\Z^d}$,
    \begin{align}
    \text{ if }f(x)\in P \text{ then }\quad \forall t\ge 0,~A_t^{\Delta_x} \subset A^{f}_t. \label{pire configuration de départ}
    \end{align}
\end{itemize}
\end{lemma}
\begin{proof} The first item comes from the additivity property: 
\[
\text{if }f\leq g, \text{ then }\xi_t^f\le \xi_t^f \vee \xi_t^g = \xi_t^{f\vee g}\le \xi_t^g.
\]
For the second item, by definition of $\Delta_x$ and since $f(x)\in P$, then $\Delta_x\leq f$. So, by the attractivity property \eqref{Propriété d'attractivité},  we have that:
\[
\forall t\ge 0,\quad\xi_t^{x} \le \xi^{f}_t.
\]
Now, let $y\in A_t^{x}$. Since $\xi_t^{x}(y) \le \xi^{f}_t(y)$ and $P$ is an upper set, then $y\in A^{f}_t$.
\end{proof}

\subsubsection{Probability space and translations} \label{Section espace de proba}

We denote by $\mathbb E^d$ the set of oriented edges between nearest neighbors of $\Z^d$. Let $((\xi_t),P)\in \mathcal C$, $(\mathcal N_{e})_{e\in \mathbb E^d}$ and $(\mathcal N_{x})_{x\in \Z^{d}}$ be mutually independent Poisson point processes on $\R^+\times \mathcal{M}$ such that the process $(\xi_t)$ can be constructed from these families. We endow $\R^+\times \mathcal M$ with the $\sigma$-algebra $\mathfrak{B}(\R^+)\otimes \mathcal{P}(\mathcal{M})$, and we denote by $M$ the set of locally finite counting measures $m=\sum_{i=0}^{+\infty} \delta_{(t_i,m_i)}$ on $\R^+\times \mathcal M$. We endow $M$ with the $\sigma$-algebra $\mathcal{F}_M$ generated by the maps $m\mapsto m(B)$, with $B\in \mathfrak{B}(\R^+)$. We define the measurable space $(\Omega,\mathcal{F})$ by 
$$(\Omega,\mathcal{F})=(M^{\E^d}\times M^{\Z^d},\mathcal{F}_M^{\otimes\E^d}\otimes\mathcal{F}_M^{\otimes \Z^d}).$$
We define the probability measure $\P=\left(\bigotimes_{e\in \mathbb E^d}\mathcal{N}_{e}\right)\otimes\left(\bigotimes_{x\in\Z^d}\mathcal{N}_{x}\right)$ on this measurable space.
For any $t\geq 0$, we define the time translation operator $\theta_t$ on a locally finite counting measure $m=\sum_{i=1}^\infty \delta_{t_i}$ on $\R_+$ by setting 
\[ \theta_t m=\sum_{i=1}^\infty \1_{\{t_i\ge t\}}\delta_{t_i-t}.\]
It induces an operator on $\Omega$, still denoted $\theta_t$: we call it \emph{time translation}. It is defined as follows:
$$\theta_t((\omega_{e})_{e\in \mathbb E^d},(\omega_{x})_{x\in \Z^d})=((\theta_t\omega_{e})_{e\in \mathbb E^d},(\theta_t\omega_{x})_{x\in\Z^d}).$$

The Poisson point processes being translation invariant, the probability measure $\mathbb{P}$ is stationary under the action of $\theta_t$. For all $z\in \Z^d$, we also define a \emph{spatial translation} $T_z$ on $\Omega$ by:
$$T_z((\omega_{e})_{e\in \mathbb E^d},(\omega_{x})_{x\in \Z^d})=((\omega_{(x+z,y+z)})_{e=(x,y)\in \mathbb E^d},(\omega_{x+z})_{x\in\Z^d}).$$

Note that if $(\xi_t)$ is in class $\mathcal{C}$, then for all $z\in \Z^d$, the processes $(\xi_t\circ T^{-z})$ and $(\xi_t)$ have the same law.

\subsubsection{Extinction time and hitting time}
We define, for all $f\in S^{\Z^d}$ and $y\in \Z^d$:
\begin{align*}
  \tau^f =\inf\{t>0,~A_t^f=\emptyset\} \quad \text{and}\quad t^f(y) =\inf\{t>0,~y\in A_t^f\},
\end{align*}
and similarly to the notations for the process $(\xi_t)$, we denote by $\tau^x$ (resp. $t^x(y)$) the quantity $\tau^{\Delta_x}$ (resp. $t^{\Delta_x}(y)$), and $\tau=\tau^0$ (resp. $t(y)=t^0(y)$). Starting from $f$, $\tau^f$ is the \emph{extinction time} of the property, that is, the first time when no one satisfies $P$; by the Definition~\ref{classeC} of the class $\mathcal{C}$ models, we have that for all $t> \tau^f$, $A_t^f=\emptyset$. The quantity $t^f(y)$ is the \emph{hitting time} of $y$, that is, the first time when $y$ satisfies $P$, starting from $f$. Finally, we denote by $\overline \P_x$ the probability measure $\P(.|\tau^x=+\infty)$, that is the probability measure conditioned on the survival of the property, starting from $\Delta_x$, and $\overline{\P}=\overline \P_0$.







\subsubsection{Linear growth models}\label{section_LGM} We work on a subclass $\mathcal C_L$ of $\mathcal{C}$ which contains processes for which the set of sites satisfying the property grows linearly with time, in the sense of the following definition.

\begin{definition}\label{classeL} We say that $((\xi_t),P)$ belongs to the class $\mathcal C_L$ if, for any finite support initial configuration $f$, it satisfies the following conditions:
\begin{enumerate}
\item we have $\P(\tau^f=+\infty)>0,$
\item there exist $C_1,C_2,M_1,M_2>0$ such that for all $t\in \R_+$ and $x\in\Z^d$,
\begin{align*}
\P\left(\exists y\in\Z^d: t^f(y)\leq t\text{ and }\|y\|_1\ge M_1 t\right)&\leq C_1\exp(-C_2t), \tag{AML} \\
\P\left(t<\tau^f<\infty\right)&\leq C_1\exp(-C_2 t),\tag{SC} \\
\P\left(t^f(x)\geq M_2\|x\|_1+t,\tau^f=\infty\right)&\leq C_1\exp(-C_2t).\tag{ALL} 
 \end{align*}
\end{enumerate}
\end{definition}

Condition~(2) means that the growth of the set of points satisfying $P$ is at most linear (AML), at least linear (ALL) and, if the extinction time $\tau$ is finite, then it is small (SC for small clusters, in analogy with percolation vocabulary). A model in $\mathcal{C}\cap \mathcal{C}_L$ is called a \emph{random linear growth model}. 

For those familiar with this literature, we briefly review here how these conditions are obtained. The first step is to create a Bezuidenhout-Grimmett block construction~\cite{BG90}, i.e. find an event in a finite space-time box whose probability large enough ensures the positive probability of survival of the process by coupling it with an oriented percolation. This construction has been implemented in~\cite{CPA},~\cite{steifwar}, \cite{durschi} and~\cite{DOP}. The second step is to perform a restart procedure: from any initial configuration, we wait until we have a good configuration before starting the block construction, then we use the known results on the oriented percolation below to obtain the desired estimates (we restart the percolation if it dies out before the top process). This procedure has been implemented in~\cite{CPA}. It is generalized in \cite{DOP} and can be applied to prove that the models from~\cite{steifwar} and \cite{durschi} satisfy the (SC) and (ALL) conditions.

\subsection{Results and examples}

We will prove the following general asymptotic shape theorem.

\begin{theorem}\label{tfa} Let $(\xi_t)$ be a Markov process taking values in $S^{\Z^d}$ and $P$ an associated property of interest. If $((\xi_t),P)$ belongs to classes $\mathcal C$ and $\mathcal C_L$, then there exists a norm $\mu$ on $\R^d$ such that for all $\varepsilon >0$, $\P(.| \tau=+\infty)$ almost surely, for all $t$ large enough, 
\[(1-\varepsilon)B_{\mu}\subset\frac{\tilde{H}_t}t\subset(1+\varepsilon)B_{\mu},\]
with $\tilde{H}_t=H_t+[0,1)^d$, and $B_\mu$ the unit ball associated to $\mu$.
\end{theorem}

The proof of this theorem is given in Section~\ref{section_thm}. In what follows, we give examples of models to which this theorem applies. This asymptotic shape theorem is already proved for examples~\ref{ex_cp} and~\ref{ex_cpa}. In each statement mentioned below, we denote by $B_\mu$ the unit ball associated to a norm $\mu$ and $\Pbarre$ the probability conditioned to survival of the property $P$ starting from the configuration $\Delta_0$.

\subsubsection{Classical contact process}\label{ex_cp} The generator of the contact process with infection rate $\lambda>0$ and recovery rate $1$ can be written in the form \eqref{Forme du générateur}, with $S=\{0,1\}$ and
\begin{align*}
c(x,\xi,0)&=\1_{\xi(x)=1}\\
c(x,\xi,1)&=\1_{\xi(x)=0}\times \lambda\sum_{y\sim x}\1_{\xi(y)=1}.
\end{align*}

As detailed in the introduction, the asymptotic shape theorem has already been proved for $P=\{1\}$, that is the property 'to be infected', by Bezuidenhout, Durrett, Griffeath and Grimmett, see \cite{DG82} and \cite{BG90}.

The next examples are extensions of the classical contact process from the literature which we define using the notation from our setting. They belong to $\mathcal C$ by construction, and with quite some effort, several authors showed that they are $\mathcal C_L$.

\subsubsection{Contact process with aging (CPA)}\label{ex_cpa} Here, we interpret the contact process dynamics as birth/death instead of infection/recovery. \cite{CPA} introduced the contact process with aging where particles have an integer 'age' affecting their birth rate so $S=\{0,1,...\}$. The parameters of the model are the aging rate $\gamma>0$ and the sequence of infection rates $(\lambda_i)_{i\geq 1}$ which is non-decreasing and limited when $i$ goes to infinity. The jump intensities in the generator are
\begin{align*}
c(x,\xi,0)&=\1_{\xi(x)>0},\\
c(x,\xi,1)&=\1_{\xi(x)=0}\times \sum_{i\geq 1}\lambda_i N_i(x,\xi), \\
\text{and } c(x,\xi,k)&=\gamma\1_{\xi(x)=k-1},
\end{align*}
where $k\geq 2$ and $N_i(x,\xi)$ is the number of neighbors of $x$ with age $i$ in $\xi$. Let $P=\{1,...\}$ be the property 'to be alive'. The process $(A_t)$ represents the set of points alive regardless of their age. Here, $\Delta_0$ is the configuration where all particles are dead $(0)$ except the origin which is alive with age $1$. From the results obtained in~\cite{CPA}, $((\xi_t),P)$ belongs to $\mathcal C \cap \mathcal C_L$ and verifies the asymptotic shape theorem in the supercritical phase. This work has inspired our proof of an asymptotic shape theorem for a general class of random linear growth models. We deduce from this result an asymptotic shape theorem for Krone's model~\cite{krone}, which is a particular case of the contact process with aging.

\subsubsection{Richardson's model with stirring (RMS) and contact process with stirring (CPS)}\label{CPS} Richardson's model is a contact process without healing. The stirring dynamic corresponds to an exchange of the states of two sites. For both RMS and CPS, we have $S=\{0,1\}$ and $P=\{1\}$ the property 'to be infected'. Here, the initial configuration $\delta_{0}$ is the classical one where all sites are healthy except the infected origin. For these two models, stirring represents the movements of the infected particles. The generators of RMS and CPS can be written in the form \eqref{Forme du générateur}. For RMS with infection rate $\lambda>0$ and stirring rate $1$, the jump intensities are:
\begin{align*}
c(x,\xi,1)&=\1_{\xi(x)=0}\times \lambda\sum_{y\sim x}\1_{\xi(y)=1},\\
c(x,y,\xi,1,0)&=\1_{\{\xi(x)=1,\xi(y)=0\}}.
\end{align*}
For CPS with infection rate $\lambda>0$, healing rate $1$ and stirring rate $N>0$, the jump intensities are:
\begin{align*}
c(x,\xi,1)&=\1_{\xi(x)=0}\times \lambda\sum_{y\sim x}\1_{\xi(y)=1},\\
c(x,y,\xi,1,0)&=N\times \1_{\{\xi(x)=1,\xi(y)=0\}},\\
c(x,\xi,0)&=\1_{\xi(x)=1}.
\end{align*}

Marchand, Marcovici and Siest~\cite{MMS} prove that both RMS and CPS are in the class $\mathcal C \cap \mathcal C_L$ for large enough infection rates, therefore Theorem~\ref{tfa} can be applied to obtain an asymptotic shape theorem in these cases. We denote by $\lambda_c(d)$ the critical infection parameter for the classical contact process with healing rate $1$ in dimension $d$.

\begin{corollary}[Asymptotic shape theorem for RMS and CPS]\label{tfa_RMS_CPS}

\begin{enumerate}
    \item Let ($\xi_t$) be a RMS with infection rate $\lambda>2d\lambda_c(d)$ and stirring rate $1$. There exists a norm $\mu$ on $\R^d$ such that for all $\varepsilon>0$, $\P$-almost surely, for all $t$ large enough, 
 \[(1-\varepsilon)B_\mu\subset\frac{\tilde{H}_t}{t}\subset(1+\varepsilon
)B_\mu,\]
where $\tilde{H}_t=H_t+[0,1)^d$, with $H_t$ the set of particles infected before time $t$. 
    \item Let ($\xi_t$) be a CPS of stirring rate $\nu$, infection rate $\lambda>(2d\nu+1)\lambda_c(d)$ and healing rate $1$. There exists a norm $\mu$ on $\R^d$ such that for all $\varepsilon>0$, $\P(.|\tau=+\infty)$-almost surely, for all $t$ large enough, 
 \[(1-\varepsilon)B_\mu\subset\frac{\tilde{H}_t}{t}\subset(1+\varepsilon
)B_\mu,\]
where $\tilde{H}_t=H_t+[0,1)^d$, with $H_t$ the set of particles infected before time $t$. 
\end{enumerate}
\end{corollary}


\subsubsection{Others examples} In~\cite{Foxall2016}, Foxall gives two other examples belonging to class $\mathcal C$: the \textit{bipartite infection model} and the \textit{household model}. The author asserts that the Bezuidenhout-Grimmett construction can easily be extended to these models, so that it would be possible to use~\cite{DOP} techniques to deduce the three conditions for belonging to $\mathcal C_L$.

\subsection{Discussion on non-additive models}

Additivity and attractivity are key ingredients in our proofs, in particular for the construction and properties of the essential hitting time. We believe that the additivity assumption is not necessary, but we need it in the proof presented in this paper. Indeed, this assumption allows us to use a percolative path structure; more precisely, it allows us to restart the process from a point that we know to satisfy the property $P$, without considering the environment, and, conversely, to consider any point satisfying the property as descending from a configuration with a single point satisfying the property. The attractivity hypothesis, meanwhile, is fundamental to all the tools and results in the paper.

Nevertheless, there are non-attractive models for which an asymptotic shape theorem is conjectured but for which, to our knowledge, proofs only hold in dimension $1$. For instance, kinetically constrained models which are interacting particle systems on $\{0,1\}^{\Z^d}$ where each site can update its state (at rate $q$ in state $1$, at rate $1-q$ in state $0$) if it satisfies a local condition (for FA1-f the condition is "to have at least one neighbor in state 1"; for East model in dimension 1, the condition is "to have its right neighbor in state 1"). These models are not attractive. Blondel~\cite{front_east}, for East model, and then Blondel, Deshayes and Toninelli~\cite{BDT}, for FA1-f with large $q$, have shown linear edge growth in dimension 1, but the question remains open in higher dimensions. Hartarsky and Toninelli~\cite{HT24} have recently proved that kinetically constrained models exhibit at least linear growth. 

Velasco~\cite{Vel} introduced two new extensions of the contact process. In both models, $S=\{-1,0,1\}$, where $-1$ state represents sterile particle (unable to give birth):
\begin{itemize}
    \item the inherited sterility process (IS), where $1$'s give birth to $1$ at rate $\lambda$ with probability $p$, and to $-1$ with probability $1-p$.
    \item the spontaneous sterility process (Spont) where $1$'s give birth to $1$ at rate $\lambda$ (as in the classical CP), but where $-1$'s can appear spontaneously in the place of 0's (without neighborhood condition).
\end{itemize}
Spont process, which can be seen as a continuous version of the Garet-Marchand bacteria model~\cite{DOP}, is attractive but non additive. For these models, (ALL), (AML) and (SC) have been established in~\cite{DOP} but the lack of additivity prevents us from applying our result to deduce an asymptotic shape theorem. The IS is also not attractive so our techniques deeply fail. In a forthcoming paper~\cite{avecisa}, the authors prove a shape theorem for these two processes in dimension 1.

Another extension of the contact process that is relevant to study is the contact process in a randomly evolving environment; it was introduced by Broman~\cite{broman} and Steif and Warfheimer~\cite{steifwar} made the block construction that could prove (ALL) and (SC). Unfortunately, the lack of additivity once again prevents us from applying our result. A recent preprint~\cite{seiler_reitmeier} addresses this problem. 

From now on, $((\xi_t),P)$ belongs to $\mathcal C \cap \mathcal C_L$. In Section~\ref{section_essential}, we introduce the essential hitting time $\sigma$, a quantity with good properties for the induced dynamical system, and show that the difference between $\sigma$ and the hitting time $t$ is not too large. In Section~\ref{section_thm}, we prove Theorem~\ref{tfa}, thanks to a subadditive ergodic theorem applied to the essential hitting time.

\section{Essential hitting time}\label{section_essential}

\subsection{Definition and properties}

With non-permanent models like extensions of the contact process, the hitting times can be infinite (because extinction is possible), and survival conditioning can cause a loss of independence, stationarity and even subadditivity properties required by Kingman's theory. We are inspired by the construction of Garet and Marchand~\cite{GMPCRE}, for the contact process in random environment, to build the \emph{essential hitting time} for a random linear growth model $((\xi_t),P)\in\mathcal C \cap \mathcal C_L$. In the following, we will use the different properties of $(\xi_t)$ detailed in Definitions~\ref{classeC} and~\ref{classeL}.

We set $u_0(x)=v_0(x)=0$ and inductively define two sequences of stopping times $(u_n(x))_n$ and $(v_n(x))_n$ as follows.
\begin{itemize}
\item Suppose that $v_k(x)$ is defined. We set 
\begin{align*}
u_{k+1}(x)&=\inf\{t\geq v_k(x):x\in A_t\}.
\end{align*} 
If $v_k(x)$ is finite, then $u_{k+1}(x)$ is the first time after $v_k(x)$ when the site $x$ satisfies $P$; otherwise, $u_{k+1}(x)=+\infty$.
\item Suppose that $u_k(x)$ is defined, with $k\geq 1$. We set $v_k(x)=u_k(x)+\tau^x\circ \theta_{u_k(x)}$. If $u_k(x)$ is finite, then the time $\tau^x\circ\theta_{u_k(x)}$ is the (possibly infinite) extinction time of $P$ starting at time $u_k(x)$ from the configuration $\Delta_x$; otherwise, $v_k(x)=+\infty$.
\end{itemize}

We have $u_0(x)=v_0(x)\leq u_1(x)\leq v_1(x) \leq \ldots \leq u_i(x) \leq v_i(x) \ldots$ 
We then define $K(x)$ to be the first step when $v_k(x)$ or $u_{k+1}(x)$ becomes infinite: 
\[
K(x)=\min\{k\geq 0: v_k(x)=+\infty\text{ or }u_{k+1}(x)=+\infty\}. \]
The following lemma says that if the property survives in the entire process $(\xi_t)$, then the procedure stops after a finite number of steps almost surely, and at a time $u_K(x)$ such that the site $x$ verifies the property $P$, and spreads it forever. 
\begin{lemma} \label{Temps d'atteinte essentiel bien défini}
    Let $\rho:=\P(\tau=+\infty)$. We have:
    \begin{enumerate}
        \item\label{K} $\forall x\in \Z^d,~\forall k\in \N,~\P(K(x)>k)\le (1-\rho)^k$.
        \item Almost surely, for all $x\in \Z^d$,
        $$(K(x)=k \text{ and }\tau=+\infty) \iff (u_k(x)<+\infty \text{ and } v_k(x)=+\infty).$$
    \end{enumerate}
\end{lemma}
\begin{proof}
\underline{First item.} Since $u_{k+1}(x)$ is a stopping time, we can apply the strong Markov property at time $u_{k+1}(x)$. Therefore we have:
\begin{align*}
    \P(K(x)> k+1)&=\P(u_{k+2}(x)<+\infty)\\
    &\le \P(u_{k+1}(x)<+\infty,v_{k+1}<+\infty)\\
    &\le \P(u_{k+1}(x)<+\infty,\tau^x\circ \theta_{u_{k+1}(x)}<+\infty)\\
    &\le \P(u_{k+1}(x)<+\infty)\P(\tau^x<+\infty)\\
    &\le \P(K(x)>k)(1-\rho),
\end{align*}
the last line is due to translation invariance.

\noindent\underline{Second item.} Suppose that $K(x)=k$ and $v_k(x)<+\infty$. Since $v_k(x)$ is a stopping time, we can apply the strong Markov property at time $v_k(x)$. Therefore, denoting by $\eta$ the configuration $\xi_{v_k(x)}$, we obtain:
\begin{align*}
    &\P(\tau=+\infty,v_k(x)<+\infty,u_{k+1}=+\infty|\mathcal{F}_{v_k(x)})\\
    &=\1_{\{v_k(x)<+\infty\}}\P(\tau^\eta=+\infty,t^\eta(x)=+\infty)= 0.
\end{align*}
which is impossible because the growth is at least linear (ALL), so in particular the process survives strongly (that is each point is reached). 
Now suppose that $u_k(x)<+\infty \text{ and } v_k(x)=+\infty$. By definition, we have $K(x)=k$ and $\tau^x\circ \theta_{u_k(x)}=+\infty$. Since the model is attractive, and since $x$ satisfies the property P at time $u_k(x)$, then we have $\tau=+\infty$, and the second implication holds.
\end{proof}
Now we can define our essential hitting time $\sigma(x)$, and a translation related to it. See Figure \ref{fig:constructionSigma} for an illustration of the construction.
\begin{definition}\label{defsigma} For $x\in\Z^d$, we call \emph{essential hitting time of $x$ by the property $P$} the quantity $\sigma(x)=u_{K(x)}$. We define the operator $\tilde{\theta}_x$ on $\Omega$ by setting:
\begin{equation*}
\tilde{\theta}_x=\begin{cases}
  T_x\circ\theta_{\sigma(x)} & \text{if } \sigma(x)<+\infty,\\
  T_x &\text{otherwise}.
\end{cases}
\end{equation*}
\end{definition}

\begin{figure}
    \centering
    \begin{tikzpicture}[xscale=0.6,yscale=0.4]

 \draw [->] (-2,-1)--(10,-1);
 \draw (10,0) node[below right] {$\Z^d$};
 \draw [->,gray] (-1.5,-2)--(-1.5,5);
 \draw[gray] (-1.5,5) node[above] {time};
  \draw (0,-1) node[below] {$0$};
  
  \draw[dashed] (7,-1)--(7,16);

\draw[->,>=latex] (0,-1) to[out=90,in=-135] (2,0);
\draw[->,>=latex] (2,0) to[out=45,in=-90] (0.5,1.1);
\draw[->,>=latex] (0.5,1.1) to[out=90,in=-120] (4,3.2);

\draw[->,>=latex](4,3.2) to[out=120,in=-90] (2,6);
\draw[->,>=latex](2,6) to[out=90,in=-110] (3.7,8);
\draw[->,>=latex](3.7,8) to[out=70,in=-110] (5,10);
\draw[->,>=latex](5,10) to[out=60,in=-100] (7,11);
\draw[->,>=latex](5,10) to[out=80,in=-120] (7,12.5);
\draw[-](7,11) to[out=80,in=-90] (9,12.5);
\draw[thick] (9,12.5)--++(0.1,0.1);\draw[thick] (9,12.5)--++(0.1,-0.1);\draw[thick] (9,12.5)--++(-0.1,-0.1);\draw[thick] (9,12.5)--++(-0.1,+0.1);

\draw[-](7,11) to[out=100,in=-90] (5,13.5);
\draw[thick] (5,13.5)--++(0.1,0.1);\draw[thick] (5,13.5)--++(0.1,-0.1);\draw[thick] (5,13.5)--++(-0.1,-0.1);\draw[thick] (5,13.5)--++(-0.1,+0.1);

\node[draw,circle,fill=red,opacity=0.4] at (7,12.5){};
\draw[dotted] (7,11)--(12,11);
\draw (12,11) node [right]{$u_2(x)$};
\draw[dotted] (5,13.5)--(12,13.5);
\draw (12,13.5) node [right]{$v_2(x)$};
\draw[->,>=latex](7,12.5) to[out=50,in=-90] (9.5,14);
\draw[->,>=latex](9.5,14) to[out=90,in=-50] (7,15.5);
\draw[dotted] (7,15.5)--(12,15.5);
\draw (12,16) node [right]{$u_3(x)=\sigma(x)$};
\draw (12,15) node [right]{$v_3(x)=\infty,~K(x)=3$};
\draw[->,>=latex](9.5,14) to[out=90,in=-120] (11.5,16.5);
\draw[->,>=latex](7,15.5) to[out=130,in=-50] (5,17);
\draw (5,17) node [above left]{$\infty$};
\draw (11.5,16.5) node [above right]{$\infty$};

\draw[->,>=latex](4,3.2) to[out=60,in=-90] (7,5);
\draw[dotted] (7,5)--(12,5);
\draw (12,5) node [right]{$u_1(x)=t(x)$};
\draw[-] (7,5) to[out=90,in=-90] (10,7);

\draw[thick] (10,7)--++(0.1,0.1);
\draw[thick] (10,7)--++(0.1,-0.1);
\draw[thick] (10,7)--++(-0.1,-0.1);
\draw[thick] (10,7)--++(-0.1,+0.1);

\draw[-] (7,5) to[out=90,in=-80] (8,8);

\draw[thick] (8,8)--++(0.1,0.1);\draw[thick] (8,8)--++(0.1,-0.1);\draw[thick] (8,8)--++(-0.1,-0.1);\draw[thick] (8,8)--++(-0.1,+0.1);

\draw[-] (7,5) to[out=100,in=-90] (5,6.5);
\draw[thick] (5,6.5)--++(0.1,0.1);\draw[thick] (5,6.5)--++(0.1,-0.1);\draw[thick] (5,6.5)--++(-0.1,-0.1);\draw[thick] (5,6.5)--++(-0.1,+0.1);

\draw[dotted] (8,8)--(12,8);
\draw (12,8) node [right]{$v_1(x)$};

\draw[<->,dotted] (0,5)--(0,15.5) node[midway, left ]{$\sigma(x)-t(x)$}; 
\draw(7,-1) node[below]{$x$};

\end{tikzpicture}
    \caption{Construction of the essential hitting time. \small At time $u_1(x)$, the site $x$ verifies the property for the first time. At time $v_1(x)$, the process starting from $\Delta_x$ at time $u_1(x)$ does not contain any sites satisfying the property: since it is an absorbing state, the property disappears forever in this process. However, the initial process spreads the property to $x$ at time $u_2(x)$, so the procedure restart at that time. At time $v_2(x)$, the process starting from $\Delta_x$ at time $u_2(x)$ does not contain any sites satisfying the property. Note that instead of pursuing the procedure at the red circle, we start at the first time \emph{after the extinction time $v_2(x)$} when $x$ satisfies the property. If we do not proceed as such, then the process obtained does not have the independence properties we want. Finally, the process starting from $\Delta_x$ at time $u_3(x)$ spreads the property forever, therefore the procedure stops.}
    \label{fig:constructionSigma}
\end{figure}

Note that, by the symmetry property of the transition rates of our processes \eqref{symmetry property}, $\sigma(x)$ and $\sigma(-x)$ have the same law. In the following lemma, we state that the essential hitting time $\sigma(x)$ and its associated translation $\tilde{\theta}_x$ verify the independence and invariance properties we need to prove our asymptotic shape theorem. We recall that we denote by $\overline \P$ the probability measure $\P(.|\tau=+\infty)$, the probability conditioned to the survival of the property, starting from $\Delta_0$.

\begin{lemma} \label{invariancePbarre}
Let $x,y\in \Z^d$, with $x \neq 0_{\Z^d}$.
\begin{enumerate}
\item[(I)] For all $A$ in the $\sigma$-algebra generated by $\sigma(x)$, and $B\in \mathcal F$, we have: 
\[
\overline{\mathbb{P}}(A \cap \tilde{\theta}_x^{-1}(B))=\overline{\mathbb{P}}(A)\overline{\mathbb{P}}(B).
\]
\item[(II)] For all $A\in \mathcal{F}$, we have $\overline{\mathbb{P}}(\tilde{\theta}_x^{-1}(A))=\overline{\mathbb{P}}(A).$
\item[(III)] Under $\overline{\mathbb{P}}$, $\sigma(y)\circ\tilde{\theta}_x$ is independent from $\sigma(x)$, and it has the same law as $\sigma(y)$.
\item[(IV)] Under $\overline{\mathbb{P}}$, the random variables $(\sigma(x) \circ(\tilde\theta_{x})^j)_{j\ge 0}$ are independent and identically distributed.
\end{enumerate}
\end{lemma}

\begin{proof}
For the first point, the proof is very similar to the one of Lemma 8 of~\cite{GMPCRE}. Simply show that, for all $k\in\N^*$, we have
\[
\overline{\mathbb{P}}(A \cap \tilde{\theta}_x^{-1}(B)\cap \{K(x)=k\})=\overline{\mathbb{P}}(A\cap \{K(x)=k\})\overline{\mathbb{P}}(B).
\]
Let $A'\in \R$ be a Borel set such that $A=\{\sigma(x)\in A'\}$. We have:
\begin{align}
    &\P(\{\tau=+\infty\} \cap A \cap \tilde{\theta}_x^{-1}(B)\cap \{K(x)=k\})\notag \\
    &=\P(\tau=+\infty, \sigma(x)\in A', (T_x\circ\theta_{\sigma(x)})^{-1}(B),u_k(x)<+\infty,v_k(x)=+\infty) \label{v_k infini et pas u_k}\\
    &=\P(u_k(x)\in A', (T_x\circ\theta_{u_k(x)})^{-1}(B),u_k(x)<+\infty,\tau^x\circ \theta_{u_k(x)}=+\infty) \label{si le sous-processus survit, alors le gros aussi}\\
    &=\P(u_k(x)\in A',u_k(x)<+\infty)\P(\{\tau^x=+\infty\}\cap T_x^{-1}(B)) \label{Propriété de Markov forte}\\
    &=\P(u_k(x)\in A',u_k(x)<+\infty)\P(\{\tau=+\infty\}\cap B).\notag
\end{align}
Line \eqref{v_k infini et pas u_k} comes from the second point of Lemma \ref{Temps d'atteinte essentiel bien défini}. The attractivity property implies that, for any stopping time $T$, we have
$$\{T<+\infty,\xi_T(x)\in P,\tau\circ T_x\circ \theta_T=+\infty\}\subset \{\tau=+\infty\},$$
from which we deduce line \eqref{si le sous-processus survit, alors le gros aussi}. We obtain line \eqref{Propriété de Markov forte} with the strong Markov property at time $u_{k}(x)$. Dividing by $\P(\tau=+\infty)$, we obtain:
\begin{align*}
    \overline\P(A \cap \tilde{\theta}_x^{-1}(B) \cap \{K(x)=k\})=\beta\overline\P(B),\\
\end{align*}
where $\beta$ is a constant that does not depend on $B$. Taking $B=\Omega$, we obtain:
\begin{align*}
    \beta=\overline\P(A \cap \{K(x)=k\}),
\end{align*}
and item~$(I)$ follows. Items~$(II)$, $(III)$ and $(IV)$ follow from item~$(I)$, see the proof of Corollary 9 of~\cite{GMPCRE} for details. 
\end{proof}

Item~$(2)$ of Lemma \ref{invariancePbarre} implies that, for all $x\in \Z^d$, $(\Omega,\mathcal{F},\overline{\P},\tilde{\theta}_x)$ is a dynamical system. Now, our goal is to prove that the essential hitting time is not too large. Since on $\{\tau=+\infty\}$, we have almost surely
\begin{align}
\sigma(x)=u_1(x)+\sum_ {i=1}^{K(x)-1}(v_i(x)-u_i(x))+\sum_ {i=1}^{K(x)-1}(u_{i+1}(x)-v_i(x)), \label{Temps d'atteinte essentiel découpé en deux sommes}
\end{align}
by definition of $\sigma(x)$ and by Lemma \ref{Temps d'atteinte essentiel bien défini}, then it amounts to study these two sums. The following lemma addresses the first one.

\begin{lemma} \label{Contrôle des v_i-u_i}
    There are $C_3,C_4>0$ such that, for all $x\in \Z^d$ and $t>0$, we have:
    $$\overline{\P}(\exists i<K(x),~v_i(x)-u_i(x)>t)\le C_3\exp(-C_4t).$$
\end{lemma}

\begin{proof}
    Using the strong Markov property at time $u_i(x)$, Lemma \ref{Temps d'atteinte essentiel bien défini}, translation invariance and property (SC), we have:
    \begin{align*}
        &\overline{\P}(\exists i<K(x),~v_i(x)-u_i(x)>t)\\
        &= \frac{\P\left(\bigcup_{i=1}^{+\infty}\{v_i(x)-u_i(x)>t\}\cap \{K(x)>i\}\cap \{\tau=+\infty\}\right)}{\P(\tau=+\infty)}\\
        &\le \frac{1}{\rho} \sum_{i=1}^{+\infty}\P(\{t<v_i(x)-u_{i}(x)<+\infty\}\cap \{u_{i}(x)<+\infty\}\cap \{\tau=+\infty\})\\
        &\le \frac{1}{\rho} \sum_{i=1}^{+\infty}\P\left(\theta^{-1}_{u_i(x)}(\{t<\tau^x<+\infty\})\cap \{u_{i}(x)<+\infty\}\right)\\
        &\le \frac{1}{\rho} \P(t<\tau^x<+\infty)\sum_{i=1}^{+\infty}\P(u_{i}(x)<+\infty)\\
        &\le \frac{1}{\rho} \P(t<\tau<+\infty)\sum_{i=1}^{+\infty}\P(K(x)>i-1)\\
        &\le \frac{1}{\rho^2} \P(t<\tau<+\infty)\\
        &\le C_3 \exp(-C_4t),
    \end{align*}
    with $C_3=\frac{1}{\rho^2} C_1$ and $C_4=C_2$.
\end{proof}

\subsection{Bad growth points}

The second sum in Equation~\eqref{Temps d'atteinte essentiel découpé en deux sommes} is harder to bound: the "reinfection time" $u_{i+1}(x)-v_i(x)$ depends on the configuration at time $v_i(x)$, therefore the strong Markov property will not be sufficient. We will introduce the notion of \emph{bad growth points} in order to treat this case: if there is none of them at time $v_i(x)$ in a certain box centered on $x$ with high probability, then the site $x$ will be reinfected quickly enough, which bounds the reinfection time $u_{i+1}(x)-v_i(x)$. 

Let $\nu_y$ be the counting measure giving the times of all possible changes of state of $y$. Thanks to the construction with Poisson point processes, the measure $\nu_y$ can be expressed as the law of a Poisson point process on $\R_+$. For every $x\in\Z^d$ and $t>0$, we say that $y$ has \emph{bad growth} (with respect to $x$ parameterized by $t$) if at least one of the following events occurs:
\begin{enumerate}
    \item\label{croissance trop rapide} the spread of the property starting from $y$, in some state $p\in P$, is faster than expected by (AML),
    \item $y$ in some state $p\in P$ has a finite but too long offspring,
    \item $y$ in some state $p\in P$ has an infinite offspring but takes too long to reinfect $x$,\label{survit mais réinfecte trop lentement}
    \item\label{Esurvittrop} there is not a single event occurring at $y$ in the time interval $[0,t/2]$.
\end{enumerate}
We denote by $E^y(x,t)$ the corresponding event:
\begin{align*}
E^y(x,t)&=\{\exists p \in P: H^{\Delta_y^p}_{t}\not\subset y+B_{M_1t}\} \tag{1}\\
  &~\cup\{\exists p \in P: ~t/2<\tau^{\Delta_y^p}<+\infty\}\tag{2} \\
  &~\cup\{\exists p \in P:\tau^{\Delta_y^p}=+\infty\text{ and} \inf\{s\ge2t: x\in A^{\Delta_y^p}_{s}\}> \kappa t\} \tag{3}\\
  &~\cup\{\nu_y[0, t/2]=0\}. \tag{4}
\end{align*}

where $\kappa=3M_1(1+M_2)$, $M_1$ and $M_2$ are the constants respectively given in (AML) and (ALL). For a laps time $L>0$, the measure $\nu_y$ allows us to count the number of bad growth points $(y,s)\in (x+B_{M_1t+2})\times[0,L]$:
\[
N_L(x,t)=\sum_{y\in x+B_{M_1 t+2}} \int_0^L \1_{E^y(x,t)}\circ\theta_s \,d (\nu_y +\delta_0)(s). 
\]
\begin{remark}\label{remark_N}
In $N_L(x,t)$, we count the number of space-time points $(y,s)$ involved in state change such that, at that moment, the event $E^y(x,t)\circ\theta_s$  occurs.
We add $\delta_0$ to the count so that, for all $x_0\in B_{M_1t+2}$, $(x_0,0)$ is counted. If $N_L(x,t)=0$, this means, in particular, that $\nu_{x_0}([0,t/2])\neq 0$. Thus, by iteration, each time interval of length $t/2$ included in $[0,L]$ above $x_0\in B_{M_1t+2} $ contains a potential arrival of $\nu$. 
  \end{remark}

Let us go back to our examples:

 \begin{itemize}
  \item For CPA we consider:
 \[\nu_y= \omega_y^1+\omega_y^\gamma+\sum_{e\in\E^d:y\in e}\omega_e^\infty,\] 
  where $\omega_y^1$, $\omega_y^\gamma$ and $\omega_e^\infty$ are the Poisson processes respectively giving the possible death times at $y$, the possible maturating and the possible birth times through an edge $e$ from $y$.
 \item For CPS, we consider:
  \[\nu_y= \omega_y^1+\sum_{e\in\E^d:y\in e}(\omega_e^{\lambda}+\omega_e^N),\]
  where $\omega_y^1$, $\omega_e^\lambda$ and $\omega_e^N$ are the Poisson processes respectively giving the possible death times at $y$ and the possible infection times or stirring times through an edge $e$ from $y$.
 \end{itemize}

We control the probability that a space-time box contains no bad growth point. 

\begin{lemma}
\label{controleboite}
There exist $A,B>0$ such that for all $L>0$, $x\in\Z^d$ and $t>0$ one has
\[\mathbb{P}\left(N_L(x,t)\ge1\right)\le A(1+L)\exp(-Bt).\]
\end{lemma}
\begin{proof}
Let $x \in\Z^d$, $t>0$ and $y \in x+B_{M_1 t+2}$. First, we control the probability of $E^y(x,t)$. There exist $A_1,B_1,A_2,B_2>0$ such that, for all $t>0$, one has
\begin{itemize}
 \item $\mathbb{P}(\exists p\in P: H^{\Delta^p_y}_{t}\not\subset y+B_{M_1t})\le A_1\exp(-B_1t)$ by attractivity and property (AML), 
 \item $\mathbb{P}(\exists p\in P:t<\tau^{\Delta^p_y}<\infty)\leq A_2\exp(-B_2t)$ by attractivity and property (SC),
 \item If $\tau^{\Delta^p_y}=+\infty$, thanks to the additivity property \eqref{Propriété d'additivité}, there exist $z\in A^{\Delta^p_y}_{2 t}$ and $p'\in P$ with $\tau^{\Delta^{p'}_z}\circ\theta_{2t}=+\infty$. 
Thus, the choice of $\kappa$ implies that
\begin{align*}
 & \left\{\tau^{\Delta^p_y}=+\infty, \inf\{s\ge2t: x\in A^{\Delta^p_y}_{s}\}> \kappa t\right\}\\
& \subset\{A^{\Delta^p_y}_{2 t}\not\subset y+B_{2M_1 t}\}
\cup\bigcup_{z\in y+B_{2M_1 t}}\{ t^{\Delta_z^{p'}}(x)\circ\theta_{2t}>(\kappa-2M_1) t ,\ \tau^{\Delta_z^{p'}}\circ \theta_{2t}=+\infty\}\\
& \subset \{A^y_{2t}\not\subset y+B_{2M_1t}\}\ \cup\\
&\quad \bigcup_{z\in y+B_{2M_1 t}}\left\{t^{\Delta_z^{p'}}(x) \circ\theta_{2t}>M_2\|x-z\|+M_1t-3M_2 ,\ \tau^{\Delta_z^{p'}}\circ \theta_{2t}=+\infty\right\}.
\end{align*}
Hence, with properties (AML) and (ALL),
\begin{align*}
 &\mathbb{P}(\tau^{\Delta_y^{p}}=+\infty, \inf\{s\ge2 t: x\in A^{\Delta_y^{p}}_{s}\}>\kappa t)\\
&\qquad\le A\exp(-2BM_1 t)+(1+4M_1t)^d A\exp\left(-B(M_1t-3M_2)\right).
\end{align*}
\end{itemize}
We also have to bound $\mathbb{P}\left(\nu_y([0, t/2])=0\right)$. Since $\nu_y$ has the law of a Poisson point process on $\R_+$, then there exists a constant $A_3>0$ such that $\mathbb{P}\left(\nu_y([0, t/2])=0\right)= \exp(- A_3t/2)$.

We obtain the existence of $A_4,B_4>0$ such that for all $x\in\Z^d$, $t>0$ and for all $y\in x+B_{M_1 t+2}$ we have
\begin{align}
    \mathbb{P}( E^y(x,t))\le A_4\exp(-B_4t)
    \label{expo_dec_E}
\end{align}

For $y\in x+B_{M_1t+2}$, we set $T^y_0=0$ and denote by $(T^y_n)_{n\ge1}$ the increasing sequence of times given by $\nu_y$. We have
\begin{align*}
 \mathbb{P}(N_L(x,t) \ge1) \le\E[N_L(x,t)] &= \sum_{y\in x+B_{M_1t+2}}\E\left[\int_0^L \1_{E^y(x,t)}\circ\theta_s \,d(\nu_y+\delta_0)(s)\right]\\
 &= \sum_{y\in x+B_{M_1t+2}}\E\left[\sum_{n=0}^{+\infty}\1_{\{T^y_n\le L\}}\1_{E^y(x,t)}\circ\theta_{T^y_n}\right]\\
  &= \sum_{y\in x+B_{M_1t+2}}\sum_{n=0}^{+\infty}\E\left[\1_{\{T^y_n\le L\}}\right] \mathbb{P}(E^y(x,t))\\
   &\le \sum_{y\in x+B_{M_1t+2}} \left(1+ \E[\nu_y([0,L])] \right) \mathbb{P}(E^y(x,t)) \\
   &\le (2M_1t+5)^d \left(1+LI\right)\mathbb{P}(E^y(x,t)),
\end{align*}
where I is the total intensity of the Poisson point process $\nu_y$.
Using Equation~\eqref{expo_dec_E} we obtain the desired upper bound. 
\end{proof}

If there are no such bad growth points, we can bound the reinfection time:

\begin{lemma}
\label{boite1} Let $x\in\Z^d$ and $t\ge 2$. If $L,s$ are positive integers such that $N_L(x,t)\circ\theta_{s}=0$ and $u_i(x)\in[s+t, s+L]$, then  $v_i(x)=+\infty$ or $u_{i+1}(x)-u_i(x)\le \kappa t$. 
\end{lemma}

\begin{proof}
If we have $v_i(x)=+\infty$, we are done; so suppose that $v_i(x)<+\infty$. Let $u$ be the last arrival of the Poisson point process $\nu_x$ before $u_i(x)$ (if $u_i(x)$ is an arrival, then $u=u_i(x)$). By the definition of $u_i(x)$, we have $x\in A_u$. Moreover, by Remark \ref{remark_N}, we have $u>u_i(x)-t/2$, and since $u_i(x)\in(s+t, s+L]$, we deduce that $u\in [s,s+L]$. Therefore, the space-time point $(x,u)$ is counted in $N_L(x,t)\circ\theta_{s}$, so the event $E^x(x,t)\circ\theta_u$ does not occur. We have:
\[
\tau^x\circ \theta_u=\tau^x\circ\theta_{u_i(x)}+u_i(x)-u=v_i(x)-u<+\infty.
\]
Moreover, time $u$ is an arrival of $\nu_x$; therefore, the event $E^x(x,t)\circ\theta_{u}$ does not occur. Thus, we have $\tau^x\circ\theta_{u}\le t/2$, and so $v_i(x)-u_i(x)\leq v_i(x)-u\leq t/2 $.

Now let us prove that we have $u_{i+1}(x)-u_i(x)\le \kappa t$. By additivity, there exists $x_0\in \Z^d$ such that $x_0\in A_{u_i(x)-t}$ and $x\in A_t^{\Delta_{x_0}^{p_0}}\circ \theta_{u_i(x)-t}$, with $p_0=\xi_{u_i(x)-t}(x_0)$. If $x_0\notin x+B_{M_1t+2}$, let $x_1\in \Z^d$ be a site on the boundary of the ball $x+B_{M_1t+2}$, and $t_1\in(u_i(x)-t,u_i(x))$ such that $x_1\in A_{t_1}$, $x\in A_{u_i(x)-t_1}^{\Delta_{x_1}^{p_1}}\circ \theta_{t_1}$, with $p_1=\xi_{t_1}(x_1)$, and $t_1$ is an arrival of $\nu_{x_1}$. The space-time point $(x_1,t_1)$ is counted in $N_L(x,t)\circ\theta_{s}$, so the event $E^{x_1}(x,t)\circ\theta_{t_1}$ does not occur. But $\|x-x_1\|>M_1t$, so it contradicts: $\forall p\in P,\ H^{\Delta_{x_1}^p}_t\subset x+B_{M_1 t}$.

So $x_0\in x+B_{M_1t+2}$. Since $N_L(x,t)\circ\theta_{s}=0$, by Remark~\ref{remark_N} we have: 
$$\nu_{x_0}([u_i(x)-t,u_i(x)-t/2])\ne 0.$$
Let $t_0$ be the first arrival of $\nu_{x_0}$ in $[u_i(x)-t,u_i(x)-t/2]$. We have $x\in A_{t-t_0}^{\Delta_{x_0}^{p'_0}}\circ \theta_{t_0}$, with $p_0'=\xi_{t_0}(x_0)$. We have $u_i(x)-t_0>t/2$, so $\tau^{\Delta_{x_0}^{p'_0}}\circ \theta_{t_0}>t/2$, and the non-occurrence of the event $E^{x_0}(x,t)\circ \theta_{t_0}$ implies that it is infinite and $\inf\{u\ge 2t: x\in A_u^{\Delta_{x_0}^{p'_0}}\circ \theta_{t_0}\}\le \kappa t$. Therefore, there exists $t_2\in [t_0+2t,t_0+\kappa t]$ such that $x\in A_{t_2}$. We proved earlier that $v_i(x)-u_i(x)\le t/2$, and we have:
$$t_2\ge t_0+2t\ge u_i(x)+t\ge v_i(x)+\frac t2,$$
so $t_2>v_i(x)$, thus $u_{i+1}(x)\le t_2$. Finally, we have:
$$u_{i+1}(x)-u_i(x)\le t_2-u_i(x)\le t_0+\kappa t-u_i(x)\le \kappa t.$$

See Figure~\ref{boitebad1} for an illustration of the proof. 
\end{proof}

Iterating Lemma~\ref{boite1}, we will be able to dominate $\sigma(x)$, which will be useful to prove almost sub-additivity in Proposition~\ref{subadd}.

\begin{figure}
    \centering
    \begin{tikzpicture}[xscale=0.70,yscale=0.65]
 \draw [->] (-2,0)--(15,0);
  \draw (0,0) node[below] {$0$};
  \draw[->,gray] (-1,-0.5)--(-1,19.5);

    \draw (-1,8.5) node{-} node[left]{$s\quad$};
      \draw (-1,9.5) node{-} node[left]{$s+t$};
\draw (-1,12.5) node{-} node[left]{$s+L$};   

\draw[->,>=latex] (0,0) to[out=90,in=-135] (2,1);
\draw[->,>=latex] (2,1) to[out=45,in=-90] (0.5,2.1);
\draw[->,>=latex] (0.5,2.1) to[out=90,in=-90] (4,3.7);
\draw[->,>=latex](4,3.7) to[out=120,in=-90] (2,6);
\draw[dotted,->,>=latex](2,6) to[out=90,in=-110] (3.7,8);
\draw[-](3.7,8) to[out=70,in=-120] (5.7,10); 

\begin{scope}[shift={(2,1)}]
\draw[->,>=latex,blue](3.7,9) to[out=70,in=-90] (3.7,9.5);
\draw[->,>=latex,blue](3.7,9.5) to[out=90,in=-180] (5,9.7);
\draw[blue] (5,9.7)--(5,10);
\draw[->,>=latex,blue](5,10) to[out=90,in=-100] (7,11);
\draw[->,>=latex](5,10) to[out=90,in=-120] (7,11.7);

\draw[-,green](7,11) to[out=80,in=-90] (8.2,11.7);

\draw[-,green](7,11) to[out=100,in=-90] (6,11.9);
\draw[dotted,green] (6,11.9)--(12,11.9);
\draw (12,11.9) node [right]{$v_i(x)$};


\draw[dotted] (3.7,9)--(12,9);
\draw (12,9) node [right]{$ u_i(x)-t$};
\draw[dotted] (7,11)--(12,11);
\draw (12,11) node [right]{$u_i(x)$};

\draw[dotted] (5,9.7)--(12,9.7);
\draw (12,9.7) node [right]{$t_0$};

\draw[<->] (13.5,11)--(13.5,11.9) node[midway, right]{\small{$\le t/2$}};
\draw[<->] (13.5,9.7)--(13.5,11) node[midway, right]{\small{$\ge t/2$}};
\draw[<->] (15,9.7)--(15,15.5) node[midway, sloped, below]{\small{$2t\le t_2\le \kappa t$}};

\draw[->,>=latex](7,11.7) to[out=50,in=-90] (9.5,14);
\draw[->,>=latex](9.5,14) to[out=90,in=-50] (7,15.5);
\draw[dotted] (7,15.5)--(12,15.5);
\draw (12,15.5) node [right]{$t_2\ge u_{i+1}(x)$};
\draw[->,>=latex](9.5,14) to[out=90,in=-110] (11,16.5);
\draw[->,>=latex](7,15.5) to[out=130,in=-50] (5,17);
\draw (10.6,16.4) node [above right]{$\infty$};

  \draw[dashed] (7,-1)--(7,16);
\draw(7,-1) node[below]{$x$};
  \draw[dashed] (3.7,-1)--(3.7,9);
  \draw(3.7,-1) node[below]{$x_0$};
 \fill[gray,opacity=0.4] (2.7,7.5)--(11.4,7.5)--(11.4,11.5)--(2.7,11.5)--cycle;
 \draw[<->] (7,7.5)--(11.4,7.5) node[midway, below]{$Mt+2$};
\end{scope}
\end{tikzpicture}
    \caption{\small The gray block does not contain any bad growth point, so we use the  4 events 
    in the bad growth definition, at arrival times, to control the reinfection time of $x$ from $x_0$. 
    }    
    \label{boitebad1}
\end{figure}

\subsection{Subadditivity and difference between $\sigma$ and $t$}

\begin{proposition} \label{subadd} 
There exists $A,B>0$ such that for all $x,y\in\Z^d$,
\[
\forall t>0, \quad\Pbarre\left(\sigma(x+y)-\left(\sigma(x)+\sigma(y)\circ
\tilde{\theta}_x\right)\ge t\right)\le A\exp\left(-B\sqrt{t}\right).
\]
Moreover, for $p\ge1$, there exists $M_p>0$ such that for all $x,y\in\Z^d$, 
\[\E[(\sigma(x+y)-(\sigma(x)+\sigma(y)\circ\tilde{\theta}_x))_+^p]\le M_p.\]
\end{proposition}

\begin{proof}
Let $x,y \in\Z^d$, $t>0$ and let $s=\sigma(x)+\sigma(y)\circ\tilde{\theta}_x$. We have:
\begin{align*}
 \overline{\mathbb{P}}\left( \sigma(x+y)>\sigma(x)+\sigma(y)\circ\tilde{\theta}_x+t\right)& \le  \overline{\mathbb{P}}\left( K(x+y)> \frac{\sqrt t}{\kappa}\right)\\
&+\overline{\mathbb{P}}\left(K(x+y) \le\frac{\sqrt t}{\kappa},\sigma(x+y)\ge s + t\right).
\end{align*}
Thanks to~\eqref{K} in Lemma~\ref{Temps d'atteinte essentiel bien défini}, the first right-hand term is bounded by $\frac{1}{\rho}\exp\left(\frac{\sqrt{t}}{\kappa}\ln(1-\rho)\right)$.
For the second term we will iterate Lemma~\ref{boite1} to prove that $\sigma$ is unlikely that big. Note that if $K(x+y) \le\frac{\sqrt t}{\kappa}$, then
$t \ge K(x+y) \kappa\sqrt t$ and 
\[\{N_{K(x+y) \kappa\sqrt t}(x+y,\sqrt t) \ge1\} \subset\{
N_{t}(x+y
,\sqrt t) \ge1\}.\]
So, 
\begin{align}
&\overline{\mathbb{P}}\left(K(x+y)\le\frac{\sqrt t}{\kappa},\sigma(x+y)\ge s + t\right)\notag\\
 &\le\overline{\mathbb{P}} \left( K(x+y) \le\frac{\sqrt t}{\kappa}, \sigma(x+y)\ge s + K(x+y)\kappa\sqrt t \right)\nonumber \\
 &\le \overline{\mathbb{P}}\left(N_t\left(x+y, \sqrt t\right)\circ\theta_{s} \ge1\right)\nonumber \\
&+\overline{\mathbb{P}}\left(\exists i<K(x+y): v_i(x+y)-u_i(x+y) >\sqrt t\right).\label{ici}
\end{align}

Indeed, if $N_{K(x+y)\kappa\sqrt{t}}(x+y,\sqrt{t}) =0$ and $\forall i\leq K(x+y),~ v_i(x+y)-u_i(x+y)\leq \sqrt{t} $ then there are two possibilities:
\begin{itemize}
    \item either $\forall i \leq K(x+y),~u_i(x+y)\leq s+\sqrt{t}$ and then $\sigma(x+y)\leq s+\sqrt{t}\leq s+K(x+y)\kappa\sqrt{t}$,
    \item or we define $i_0=\max\{i:u_i(x+y)\leq s+\sqrt{t}\}$ then $v_{i_0}(x+y)\leq s+2\sqrt{t}$ and $\tau^{x+y}\circ \theta_s=+\infty$ so $\inf\{r\geq 2\sqrt{t}, x+y\in A^{x+y}_r\}\circ \theta_s \leq \kappa \sqrt{t}$ and finally $u_{i_0+1}(x+y)\leq s+\kappa\sqrt{t}$ . Furthermore $u_{i_0+j}(x+y)> s+\sqrt{t}$ for all $j\geq 1$ and we can iterate Lemma~\ref{boite1} to conclude that for all $j\leq K(x+y)-i_0, u_{i_0+j}(x+y)\leq s +j\kappa\sqrt{t}$ and then $\sigma(x+y)\leq s+K(x+y)\kappa\sqrt{t}$.
\end{itemize}

Since $N_t(x+y,\sqrt{t})\circ{\theta_s}=N_t(0,\sqrt{t})\circ T_x\circ T_y\circ \theta_{\sigma(x)}\circ\theta_{\sigma(y)\circ\tilde{\theta}_x}$, we have
\[N_t\left(x+y,\sqrt t\right)\circ\theta_{s}=N_t\left(0,\sqrt{t}\right)\circ\tilde{\theta
}_y\circ\tilde{\theta}_x.\]

Thus,
$\overline{\mathbb{P}}(N_t(x+y,\sqrt t)\circ~\theta_s \ge1)= \overline{\mathbb{P}} (N_t(0,\sqrt t)\ge1)$, which is controlled by Lemma~\ref{controleboite}. And the second term of the sum~\eqref{ici} is bound by Lemma~\ref{Contrôle des v_i-u_i}.
\end{proof}

\begin{proposition}\label{difference1}
There are $A,B,\alpha>0$ such that for every $z>0$ and every $x\in\Z^d$, 
\[\overline{\mathbb{P}}\left(\sigma(x)\ge t(x)+K(x)\left(\alpha\log(1+\|x\|) +z\right)\right)\le A\exp(-Bz).\]
\end{proposition}

\begin{proof} 

We use the same method as in the proof of the previous proposition. Let us set $t=\frac{\alpha(\log(1+\|x\|) +z)}{\kappa}$, where $\alpha$ is a parameter bigger than 1 that will be adjusted later. In particular, $t\geq z/\kappa$.
\begin{align*}
\Pbarre\left(\sigma(x)\ge t(x)+K(x)\kappa t \right) 
& \le \Pbarre \left(\sigma(x)\ge K(x)\kappa t, \; N_{K(x)\kappa t}(x, t)=0 \right) +\Pbarre \left(N_{K(x)\kappa t}(x, t)\ge1 \right) \\
& \le \Pbarre(\exists i<K(x): \; v_i(x)-u_i(x) >t)+  \Pbarre(N_{K(x)\kappa t}(x, t) \ge 1) . 
\end{align*}
The first term is bounded by by Lemma~\ref{Contrôle des v_i-u_i}. For the second one, we have that:
\begin{align*}
\Pbarre \left(N_{K(x)\kappa t}(x,t)\ge 1\right)&\le   \sum_{j=1}^{+\infty}\sqrt{\Pbarre(K(x)=j)}\sqrt{\Pbarre(N_{j\kappa t}(x,t)\ge 1)}\\
&\le \sum_{j=1}^{+\infty}\sqrt{\Pbarre(K(x)=j)}\sqrt{A(1+j\kappa t)\exp(-Bt)}\\
&\le \sqrt { A(1+\kappa t) }\exp(-\frac{B}2 t) \sum_{j=1}^{+\infty}\sqrt{(1+j)\Pbarre(K(x)=j)} .
\end{align*}
Using the fact that $K(x)$ is sub-geometrical, the sum is finite. Then, take $\alpha$ large enough, there exist $A'$ and $B'$ such that the last term is bounded by $A'\exp(-B'z)$.

 \end{proof}

\begin{proposition}
    \label{difference2} $\overline{\mathbb{P}}$-almost surely, it holds that $\displaystyle\lim_{\|x\|\to+\infty} \frac{|\sigma(x)-t(x)|}{\|x\|}=0.$
\end{proposition}
 \begin{proof}
We want to bound $\overline{\E}(|\sigma(x)-t(x)|^p)$. Using Minkowski inequality, Proposition~\ref{difference1} and item~\ref{K} of Lemma~\ref{Temps d'atteinte essentiel bien défini}, we obtain 
\[ 
\overline{\E}(|\sigma(x)-t(x)|^p) \leq CP(\log(1+\|x\|))^p.
\]
So, the sequence $\left(\frac{|\sigma(x)-t(x)|}{1+\|x\|}\right)_{x\in\Z^d}$ is in $\ell^p(\Z^d)$ and tends to $0$ almost surely.
\end{proof}

\begin{corollary}\label{coro_subadd}
\begin{enumerate}[(I)]
\item\label{ALLsigma} There exist $A,B,C>0$ such that
\[\forall x \in\Z^d,\forall t>0,~ \overline{\mathbb{P}}\left(\sigma(x)\ge C\|x\|+t\right) \le A\exp(-B\sqrt{t}).\]
\item\label{momentsigma} For $p\ge1$, there exists a constant $C_p>0$ such that
\[\forall x\in\Z^d,~\overline{\E}[\sigma(x)^p] \le C_p (1+\|x\|)^{p}.\]
\item\label{uniforme} For every $\epsilon>0$, $\overline{\mathbb{P}}$-a.s., there exists $R>0$ such that
\[\forall x,y \in\Z^d,~(\|x\|\ge R\text{ and }\|x-y\|\le\epsilon\|x\|)\Longrightarrow\left(|\sigma(x)-\sigma(y)|\le C\epsilon\|x\|\right).\]
\end{enumerate}
\end{corollary}

\begin{proof} For the first item~\eqref{ALLsigma} we can write  
\begin{align*}
\overline{\mathbb{P}}\Big(\sigma(x)> (C+1)\|x\|+t\Big)  
&\le \overline{\mathbb{P}}\left((t(x)\ge C{\|x\|}+t/2 \right) +\overline{\mathbb{P}}\left( K(x)>\frac1{2\alpha} \sqrt{\|x\|+t/2}\right)\\
+&\overline{\mathbb{P}}\left(\sigma(x)> \left(C+1\right)\|x\|+t,~t(x)< C{\|x\|}+t/2,~ K(x)\le \frac1{2\alpha} \sqrt{\|x\|+t/2}\right) \\
&\le \overline{\mathbb{P}}\left(t(x)\ge C{\|x\|}+t/2 \right) +\overline{\mathbb{P}}\left( K(x)>\frac1{2\alpha} \sqrt{\|x\|+t/2}\right)\\
&~~+\overline{\mathbb{P}}\left(\sigma(x)> t(x) +K(x)\times 2\alpha \sqrt{\|x\|+t/2}\right).
\end{align*}
The first term is bounded by (ALL), the second one using Lemma~\ref{Temps d'atteinte essentiel bien défini} and the third one using Proposition~\ref{difference1} with $z=B\sqrt{t}$. The second item \eqref{momentsigma} comes from Minkowski inequality and~\eqref{ALLsigma}. Then, for the last item~\eqref{uniforme}, for $m \in\N$ and $\varepsilon>0$, we introduce
\[B_m(\epsilon)=\{\exists x,y \in\Z^d:\|x\|=m, \|x-y\|\le
\epsilon m \mbox{ and } |\sigma(x)-\sigma(y)|>C\epsilon m \}.\]
We have
\begin{align*}
\overline{\mathbb{P}}( B_m(\epsilon) )& \le \mathop{\sum_{(1-\epsilon) m\le\|x\| \le(1+\epsilon)m}}_{\|z\|\le \epsilon m} \overline{\mathbb{P}}\left(\sigma(z)\circ\tilde{\theta}_{x}+r(x,z)>C\epsilon m\right)\\
& \le \mathop{\sum_{(1-\epsilon) m\le\|x\| \le(1+\epsilon)m}}_{\|z\|\le \epsilon m} \overline{\mathbb{P}}\left(\sigma(z)>2C\epsilon m/3\right)+\overline{\mathbb{P}}\left(r(x,y-x)>C\epsilon m/3\right).
\end{align*}
Using \eqref{ALLsigma} and Proposition~\ref{subadd}, we have
\begin{align*}
\overline{\mathbb{P}}( B_m(\epsilon) )
& \le (1+2\epsilon m)^d\left(1+2(1+\epsilon) m\right)^dA\exp\left(-B\sqrt{C\epsilon m /3}\right)\\
& {} +A'\exp\left(-B'\sqrt{C'\epsilon m /3}\right)
\end{align*}
and we can conclude thanks to Borel-Cantelli lemma.
\end{proof}

\section{Proof of the asymptotic shape theorem}\label{section_thm}
In this section we will use the following almost subadditive ergodic theorem, which is a reformulation of the theorem of Kesten and Hammersley (see \cite{hamm74} and \cite{kestentrick}), adapted for the essential hitting time and the properties it verifies.

\begin{proposition}[Theorem 39 of~\cite{CPA}]
\label{metaTFA} Let $(\Omega,\mathcal{F},\mathbb{P})$ be a probability space. Let $(\sigma(x))_{x\in\Z^d}$ be random variables with finite second moments and such that, for every $x\in\Z^d$, $\sigma(x)$ and $\sigma(-x)$ have the same distribution. Let $(s(y))_{y\in\Z^d}$ and $(r(x,y))_{x,y\in\Z^d}$ be collections of random variables such that:
\begin{enumerate}
 \item[Hyp 1:]\label{Hssadd} $\forall x,y\in\Z^d,~\sigma(x+y)\leq\sigma(x)+s(x,y)+r(x,y)$
with $s(x,y)$ having the same law as $\sigma(y)$, and being independent from $\sigma(x)$,
\item[Hyp 2:]\label{Hcontrolereste} $\forall x,y\in\Z^d,~\exists C_{x,y}$ and $\alpha_{x,y}<2$ such that \\$\forall n,p, \E[r(nx,py)^2]\le C_{x,y}(n+p)^{\alpha_{x,y}}$,
\item[Hyp 3:]\label{Hau+lin} $\exists C>0$ such that $\forall x\in\Z^d,~ \mathbb{P}(\sigma(nx)>Cn\|x\|)\xrightarrow{n\rightarrow\infty} 0$,
\item[Hyp 4:]\label{Hunif} $\exists K>0$ such that $\forall \epsilon>0, \mathbb{P}-p.s ~\exists M$ such that $(\|x\|\geq M\text{ and }\|x-y\|\leq K\|x\|)\Rightarrow \|\sigma(x)-\sigma(y)\|\leq \epsilon\|x\|$,
\item[Hyp 5:]\label{Hau-lin} $\exists c>0$ such that $\forall x\in\Z^d,~ \mathbb{P}(\sigma(nx)<cn\|x\|)\xrightarrow{n\rightarrow\infty} 0$.
\end{enumerate}
Then there exists $\mu:\Z^d\to\R^+$ such that 
\[\lim_{\|x\|\rightarrow \infty} \frac{\sigma(x)-\mu(x)}{\|x\|}=0~ a.s.\]
Moreover, $\mu$ can be extended to a norm on $\R^d$ and we have the following asymptotic shape theorem: for all $\epsilon>0$, $\mathbb{P}$ almost surely, for all $t$ large enough,
\[(1-\epsilon)B_{\mu}\subset\frac{\tilde G_t}t\subset(1+\epsilon
)B_{\mu},\]
where $\tilde{G}_t=\{x\in\Z^d:\sigma(x)\le t\}+[0,1]^d$ and $B_\mu$ is the unit ball for $\mu$.
\end{proposition}

We now deduce the expected asymptotic shape theorem for the hitting time $t$:

\begin{proposition}\label{but} There exists a norm $\mu$ on $\R^{d}$ such that almost surely under $\overline{\mathbb{P}}$, \[\lim_{\|x\|\to+\infty} \frac{t(x)-\mu(x)}{\|x\|}=0,\]
and for every $\epsilon>0$, $\overline{\mathbb{P}}$-a.s., for every large $t$,
\[(1-\epsilon)B_{\mu}\subset\frac{\tilde H_t}t\subset(1+\epsilon)B_{\mu}\] where $\tilde{H}_t=\{x\in\Z^d/ t(x)\le t\}+[0,1]^{d}$ and $B_\mu$ is the unit ball for $\mu$. 
\end{proposition}
Theorem~\ref{tfa} is contained is the previous result.

\begin{proof} First, we use Proposition~\ref{metaTFA} to show that $\sigma$ satisfies an asymptotic shape theorem. We check the hypotheses of Proposition~\ref{metaTFA} using the controls of Corollary~\ref{coro_subadd}. Thanks to~\eqref{ALLsigma}, $\sigma$ has finite second moment required. 
We take $s(y)=\sigma(y)\circ\tilde{\theta}_x$ and $r(x,y)=\sigma(x+y)-(\sigma(x)+\sigma(y)\circ\tilde{\theta}_x))_+$. The hypotheses 1 and 2 are satisfied thanks to properties of $\tilde{\theta}_x$ (Lemma \ref{invariancePbarre}) and Proposition \ref{subadd}. The hypothesis~3 is the at least linear growth~\eqref{ALLsigma}. The hypothesis 4 is the control~\eqref{uniforme}. Finally, hypothesis~5 is immediately checked thanks to the at most linear growth (AML): 
$$\overline{\mathbb{P}}(\sigma(nx)< M_2n\|x\|) \le \overline{\mathbb{P}}(t(nx)< M_2n\|x\|) \le \frac{A}{\rho}\exp(-BM_2n\|x\|).$$ 
So, $\sigma$ satisfies a shape theorem.
We deduce the result for $t$ from the result for $\sigma$ thanks to Proposition \ref{difference2}. To deduce the geometric result from the analytical one, the classical proof can be done by contradiction: assuming that the geometric inclusion is not true, it is possible to construct a sequence $(x_n)$ of points such that $\mu(x_n)/t(x_n) $ does not converge to 1.
\end{proof}

\bibliographystyle{alpha}

\end{document}